\documentclass[12pt]{amsart}

\usepackage{epsf,amssymb}

\newtheorem{thm}{Theorem}[section]

\newtheorem{prop}[thm]{Proposition}

\newcommand{\Z}{\mathbb{Z}}

\newcommand{\bdry}{\partial}
\newcommand{\sa}{\rightsquigarrow}
\newcommand{\s}{\vskip.1in}
\newcommand{\n}{\noindent}

\newcommand{\be}{\begin{enumerate}}
\newcommand{\ee}{\end{enumerate}}

\begin{document}
\title{On the Gabai-Eliashberg-Thurston Theorem}

\author{Ko Honda}
\address{University of Southern California, Los Angeles, CA 90089} 
\email{khonda@math.usc.edu}
\urladdr{http://math.usc.edu/\char126 khonda}

\author{William H. Kazez}
\address{University of Georgia, Athens, GA 30602}
\email{will@math.uga.edu}
\urladdr{http://www.math.uga.edu/\char126 will}

\author{Gordana Mati\'c}
\address{University of Georgia, Athens, GA 30602}
\email{gordana@math.uga.edu}
\urladdr{http://www.math.uga.edu/\char126 gordana}

\date{November 8, 2001.}

\thanks{KH supported by NSF grant DMS-0072853; GM supported by 
NSF grant DMS-0072853; WHK supported by NSF grant DMS-0073029.}

\begin{abstract}

We present a new, completely three-dimensional proof of the fact, due to 
Gabai-Eliashberg-Thurston, that every closed, oriented, irreducible 3-manifold 
with nonzero second homology carries a universally tight contact structure.

\end{abstract}
\maketitle
The main result of this paper is a new three-dimensional proof of the following:

\begin{thm}[Gabai-Eliashberg-Thurston]  \label{GET-closed}
Let $M$ be a closed, oriented, connected, irreducible 3-manifold with 
$H_2(M,\Z)\not=0$.  Then $M$ carries a universally tight contact structure. 
\end{thm}

In \cite{HKM2}, the authors proved the following variant of the above theorem 
for manifolds with boundary:

\begin{thm}  \label{GET-boundary}
Let $(M,\gamma)$ be an oriented, compact, connected, irreducible, sutured 
3-manifold which has nonempty boundary, is taut, and has annular sutures.  Then 
$(M,\gamma)$ carries a universally tight contact structure. 
\end{thm}

Our interest in reproving these theorems is twofold.  First, if the starting 
point is a sutured manifold decomposition and the goal is to build a 
universally tight  contact structure, it should not be necessary (indeed it is 
not) to construct a  taut foliation, perturb it into a contact structure, and 
argue using symplectic  filling techniques that the resulting contact structure 
is universally tight.  Our other motivation is to use these theorems as 
guidelines in the development of a cut-and-paste theory of contact topology.  
This theory contrasts with foliation theory right from the start.  Given a tight 
contact structure, it is very easy to produce useful decompositions of the space 
(so-called {\it convex decompositions}, see \cite{HKM2}).  On the other hand, 
given two pieces of a manifold, each with a universally tight contact structure, 
it is surprisingly difficult to find gluing theorems which allow one to conclude 
that the contact structure on their union is universally tight.  In the theory 
of taut foliations, the relative difficulty levels of the two appear to be 
switched.

We have used the search for cut-and-paste proofs of the 
Gabai-Eliashberg-Thurston theorem as a method for finding new gluing theorems.  
Theorem~\ref{GET-boundary} of \cite{HKM2} was proved using gluing techniques 
(pioneered  by Colin~\cite{Co97,Co99a}).  The key gluing theorem used in the 
proof was:

\begin{thm}[Colin~\cite{Co99a}] \label{colin-gluing}
Let $(M,\xi)$ be an oriented, compact, connected, irreducible, contact 
3-manifold and $S\subset M$ an incompressible convex surface with nonempty 
Legendrian boundary and $\bdry$-parallel dividing set $\Gamma_S$.  If  
$(M\setminus S, \xi|_{M\setminus S})$ is universally tight, then $(M,\xi)$ is 
universally tight. 
\end{thm}

The condition $\bdry S\not=\emptyset$ is an important condition in the proof of 
Theorem~\ref{colin-gluing}.  Therefore, Theorem~\ref{colin-gluing} is not 
applicable when $M$ is a closed 3-manifold and the first cut in the sutured 
manifold  decomposition is along a closed surface.  Instead, we will make use of 
the main technical result of this paper, Theorem~\ref{closedgluing}, which is a 
gluing theorem along certain closed surfaces.

A predecessor to this gluing theorem is the gluing theorem along incompressible 
tori \cite{HKM2}, where we rephrased and gave a slightly different proof of a 
gluing result of Colin \cite{Co99a}.   Using it, we presented, among other 
things, a foliation-theory-free proof of the existence theorem for tight contact 
structures in the case of a closed, irreducible, toroidal manifold (without the 
assumption $H_2(M,\Z)\not=0$).  The two key ingredients of  
Theorem~\ref{GET-closed} in the toroidal case are the above-mentioned variant of 
the Gluing Theorem (Theorem~\ref{closedgluing}) along incompressible tori and a 
good understanding of universally tight contact structures on $T^2 \times 
[0,1]$. 

It suffices, for the purposes of this paper, to prove a gluing theorem for 
atoroidal manifolds, along closed convex surfaces $\Sigma$ of genus $g > 1$ that 
satisfy the {\em extremal condition} $$\langle e(\xi),\Sigma\rangle=\pm 
(2g-2),$$ where the left-hand side is the Euler class of $\xi$ evaluated on 
$\Sigma$. (Note that the condition is trivially satisfied for genus one 
surfaces.)  Tight contact structures $\xi$ on $M$ satisfying this condition are 
said to be {\em extremal along $\Sigma$} for the following reason.  The 
Bennequin inequality \cite{Be,E92} states that:  $$-(2g-2)\leq \langle 
e(\xi),\Sigma\rangle\leq 2g-2.$$ 

One of the main results of \cite{HKM3} is the classification of extremal tight 
contact structures on $\Sigma \times [0,1]$ in the case of two dividing curves 
on each boundary component.  This result, and its implications for covering 
spaces of $\Sigma \times [0,1]$, are enough to construct contact structures 
satisfying the hypotheses of the Gluing Theorem.

In Section~\ref{section:straddling}, we introduce the notion of {\it straddling} 
or {\it marking} and reformulate the classification theorem of \cite{HKM3} in a 
language more suitable for the application of the Gluing Theorem. 
Section~\ref{section:gluing} contains a proof of the Gluing Theorem.  In 
Section~\ref{section:get}, Gabai's well-groomed sutured manifold decomposition 
theory is used to construct a universally tight contact structure on the 
cut-open manifold, and then the Gluing Theorem completes the proof of 
Theorem~\ref{GET-closed}.

\s\n
We adopt the following conventions: 
\be

\item The ambient manifold $M$ is an oriented, compact $3$-manifold.

\item $\xi$ = positive contact structure which is co-oriented by a global 1-form 
$\alpha$.

\item A convex surface $\Sigma$  is either closed or compact with Legendrian 
boundary.

\item $\Gamma_\Sigma$ = dividing multicurve of a convex surface $\Sigma$.

\item $\#\Gamma_\Sigma$ = number of connected components of $\Gamma_\Sigma$.

\item $\Sigma\setminus \Gamma_\Sigma=\Sigma_+\cup \Sigma_-$, where $\Sigma_+$  
(resp.\ $\Sigma_-$) is the region where the normal orientation of $\Sigma$ is 
the same as (resp.\ opposite to) the normal orientation for $\xi$. 

\ee

\s\n
{\it Acknowledgements.}  The first author thanks the American Institute of 
Mathematics, Stanford University, and IHES for their hospitality.

\section{Straddling} \label{section:straddling}

The following proposition will be used repeatedly throughout the paper.

\begin{prop}\label{add=dig} 
Let $(M, \xi)$ be a tight contact 3-manifold with convex boundary, $\Sigma$ a 
component of $\partial M$, and $\gamma,\gamma'$ a pair of parallel disjoint 
curves in $\Gamma_{\Sigma}$. Suppose there is a bypass 
$\mathcal{B}_\alpha\subset M$ attached along an arc $\alpha\subset \Sigma$ that 
starts on $\gamma'$, crosses $\gamma$, and then ends on a third curve in 
$\Gamma_\Sigma - (\gamma \cup \gamma')$.  Let $\beta$ be a Legendrian arc in 
$\Sigma$ that starts on $\gamma$, crosses $\gamma'$, ends on a point of 
$\Gamma_\Sigma - (\gamma \cup \gamma')$, and does not intersect $\alpha$ or any 
other points of $\Gamma_\Sigma$.  Then attaching a bypass $\mathcal{B}_\beta$ to 
$M$ along $\beta$ produces a manifold contact isomorphic to the manifold 
obtained by removing a convex neighborhood of $\mathcal{B}_\alpha$. 
\end{prop}

\begin{figure}[ht]	
	{\epsfysize=1.8in\centerline{\epsfbox{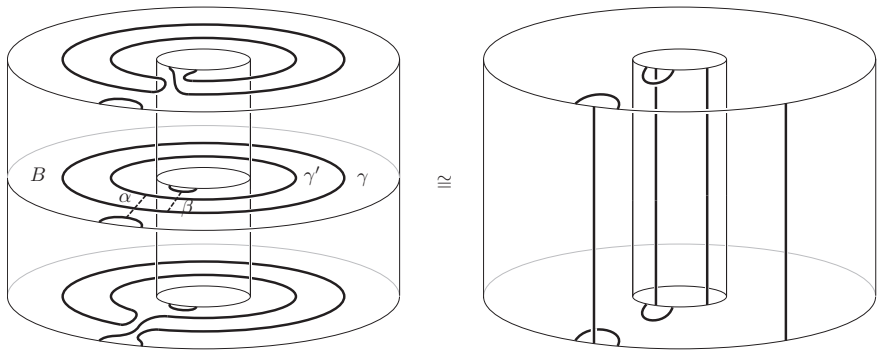}}}
	\caption{Neighborhood of $B$} 	
	\label{bypasses}
\end{figure}

\begin{proof} 
Let $B\subset \Sigma$ be a regular neighborhood of the union of $\alpha$, 
$\beta$, and the annulus in $\Sigma$ bounded by $\gamma$ and $\gamma'$ --- we 
assume $B$ is convex with Legendrian boundary.  Let $V$ be a small neighborhood 
of the union of $B$, $\mathcal{B}_\alpha$, and $\mathcal{B}_\beta$.  See 
Figure~\ref{bypasses}. Topologically, $V$ is the product of an annulus and an 
interval, i.e., a solid torus.  After the necessary edge-rounding, we see that 
the dividing set of $\partial V$ has two components, each 
of which intersects a compressing disk in a single point.  There is a unique 
tight contact structure on $V$ with this boundary condition, as can be seen by 
splitting $V$ along the compressing disk and appealing to Eliashberg's 
uniqueness theorem on a ball \cite{E92}.  It remains to verify that the contact 
structure on $V$ is indeed tight --- for this we simply remark that an explicit 
model can be found inside the unique (product) tight contact structure on $V$.  
Since the contact structure on $V$ is a product, it follows that adding a bypass 
along $\beta$ is equivalent to removing a bypass along $\alpha$. \end{proof}

Let $(M,\xi)$ be a tight contact 3-manifold with convex boundary and let 
$\Sigma$ be a connected component of $\bdry M$ of genus $g$ which satisfies 
the extremal condition $\langle e(\xi),\Sigma\rangle=- (2g-2)$.  It follows that 
$\Sigma_{-}$ has zero Euler characteristic and is a disjoint union of annuli 
$A_i$, for $i=1,\dots, n$.  Denote $\partial A_i = \gamma_i^0 \cup \gamma_i^1$, 
and call $\gamma_i^0$ and $\gamma_i^1$ a {\em parallel pair} of dividing curves. 

We say a dividing curve $\gamma_i^j$ is {\em straddled} if there is a bypass in $M$ 
with an attaching arc $\alpha_i^j$ which starts on $\gamma_i^{1-j}$, crosses 
$\gamma_i^j$, and ends on (i) a point of $\Gamma_\Sigma - (\gamma_i^0 \cup 
\gamma_i^1)$ or on (ii) $\gamma_i^0$ or $\gamma_i^1$, but only after first going 
around a nontrivial loop on $\Sigma$ which is not parallel to $\gamma_i^j$.   
Our bypasses may be degenerate, i.e., the two endpoints of the arc of attachment 
are allowed to coincide. The attaching arcs of these bypasses are called {\em 
straddling arcs}. 

A closed convex surface $\Sigma$ of genus $g>1$ which is a connected component 
of $\bdry M$ {\em admits a full marking} if the following hold:

\be
\item $\xi$ is extremal along $\Sigma$ and satisfies $\langle 
e(\xi),\Sigma\rangle=- (2g-2)$. 

\item The components of $\Sigma_-$ are pairwise nonparallel annuli.

\item There is a collection $S=\{\alpha_i\}$ of straddling arcs (called a {\em 
straddling set}) and a corresponding collection of bypasses 
$\mathbb{B}=\{\mathcal{B}_{\alpha_i}\}$ in $M$ such that:

\be
\item at least one curve in each parallel pair of $\Gamma_{\Sigma}$ is 
straddled by a bypass in $\mathbb{B}$,

\item every curve in $\Gamma_\Sigma$ is straddled by at most one bypass in 
$\mathbb{B}$,

\item if $i\not= j$, then either $\mathcal{B}_{\alpha_i}$ and 
$\mathcal{B}_{\alpha_j}$ are disjoint or intersect only at the endpoints of 
their corresponding arcs of attachment, and

\item $S$ is an {\em essential family}, i.e., $\Sigma_+ - (\cup_{i=1}^n   
\alpha_i)$ has no disk components.   (Equivalently, a thickening of 
$\Sigma_-\cup(\cup_{i=1}^n \alpha_i)$ is an incompressible subsurface of 
$\Sigma$.) 

\ee

\ee

\begin{prop}\label{covers} 
Let $(M, \xi)$ be a tight contact manifold with convex boundary, and let  
$(\widetilde M, \widetilde {\xi})$ be a finite cover.  If $\Sigma$ is a boundary 
component of $(M, \xi)$ which admits a full marking, then the preimage 
$\widetilde {\Sigma}$ of $\Sigma$ admits a full marking in $(\widetilde M, 
\widetilde {\xi})$. \end{prop}

\begin{proof}
Let $S$ be the straddling set for $\Sigma$ and $\mathbb{B}$ the corresponding 
set of bypasses.  Their preimages $\widetilde S$ and $\widetilde{\mathbb{B}}$ 
satisfy all of the axioms necessary for $\widetilde \Sigma$ to admit a full 
marking, except that there may be several arcs in {\it $\widetilde S$} which 
straddle the same dividing curve of $\widetilde \Sigma$.  Removing extra 
components of {\it $\widetilde S$} decreases the Euler characteristics of the 
complementary regions.  Thus no disk components are produced, and $\widetilde 
\Sigma$ is fully straddled. \end{proof}

\begin{prop} \label{bypass-keeps-straddled}  
Let $(M,\xi)$ be a tight contact manifold with convex boundary, $\Sigma$ a 
boundary component of $M$ which admits a full marking, and $S$ a corresponding 
straddling set.  Let $\alpha$ be an arc of $S$ which straddles the dividing 
curve $\gamma$.  Suppose there exists an arc $\alpha'\subset \Sigma$ such that:
\be
\item $S \cup\{\alpha'\}$ is an embedded, essential family of arcs, and
\item $\alpha'$ starts on $\gamma$, crosses the dividing curve $\gamma'$ 
parallel to $\gamma$, and ends on a point of $\Gamma_\Sigma - (\gamma \cup 
\gamma')$, while intersecting no other points of $\Gamma_\Sigma$.
\ee 
Then the manifold $(M', \xi')$, obtained by  attaching a bypass along $\alpha'$ 
to $(M,\xi)$, also admits a full marking along the boundary component $\Sigma'$ 
corresponding to $\Sigma$. If $(M,\xi)$ is universally tight, then so is 
$(M',\xi')$. \end{prop}

\begin{proof}
Attaching a bypass along $\alpha'$ decreases the number of dividing curves on 
$\Sigma$ by two.  The new collection of straddling arcs $S'$ is defined to be $S 
-\{\alpha\}$. To see that $S'$ is an essential family of arcs, it is necessary 
to relate the complementary region $\Omega'=\Sigma'_+ - (\cup_{\beta\in 
S'} \beta)$ to the corresponding complementary region 
$\Omega=\Sigma_+-(\cup_{\beta\in S} \beta)$.  Performing the bypass 
along $\alpha'$ is equivalent to the following two steps:
\begin{itemize}
\item First cut $\Omega$ along a subarc of $\alpha'$.  This creates no disks 
since $S \cup \{\alpha'\}$ is an essential family. 
\item Next attach a band between two disjoint intervals on the boundary of
the  modified $\Omega$. This can only decrease the Euler characteristic of the 
complementary regions.  
\end{itemize}
Finally, deleting $\alpha$ from $S$ decreases the Euler characteristic 
of the complementary regions, and hence $S'$ is an essential family.  The 
universal tightness of $(M',\xi')$ follows directly from 
Proposition~\ref{add=dig}. \end{proof}

\begin{prop} \label{straddled-unstraddled} 
Let $M$ be a closed hyperbolic manifold, $\xi$ a universally tight contact 
structure on $M$, and $\Sigma$ an incompressible, closed, convex surface in $M$ 
which satisfies the extremal condition.  Let $\Sigma^0$ be a boundary component 
of $M \setminus \Sigma$ which admits a full marking.  If $\gamma$ and $\gamma'$ 
form a parallel pair of $\Gamma_{\Sigma^0}$ and $\gamma$ is straddled, then 
$\gamma'$ is {\em unstraddled}, i.e., there is no bypass in $(M\setminus 
\Sigma,\xi)$ straddling $\gamma'$. 
\end{prop}

\begin{proof} 
We argue by contradiction.  Let $S$ be the straddling set for $\Sigma^0$, 
$\mathbb{B}$ be the corresponding set of bypasses, and $\alpha$ be the arc in 
$S$ which straddles $\gamma$.  Suppose there exists a straddling arc $\beta$ for 
$\gamma'$.  We claim there exists a finite cover $\pi:(\widetilde M, 
\widetilde{\xi}) \to (M, \xi)$ in which components $\widetilde{\alpha}$ and 
$\widetilde{\beta}$ of preimages of $\alpha$ and $\beta$ straddle a parallel 
pair of dividing curves $\widetilde{\gamma}$ and $\widetilde \gamma'$, and are 
disjoint.  To prove the claim we use the following theorem of Allman-Hamilton 
\cite{AH}.

\begin{thm}[Abelian subgroup separability]\label{AH} 
Let $M$ be a hyperbolic 3-manifold.  Then abelian subgroups $H$ of  
$\pi_1(M,*)$ are separable, i.e., for any $g \in \pi_1(M,*) - H$ there 
exists a finite index subgroup $K \supset H$ which does not contain $g$. 
\end{thm}

\n
Let $\sigma$ be the core of the annulus $A$ bounded by $\gamma$ and $\gamma'$.  
If $\alpha$ intersects only two dividing curves, let $\overline\alpha$ be a 
closed loop formed by the union of $\alpha$ and a subarc of $\sigma$.  From the 
definition of a straddling arc, $\overline\alpha$ is not in the subgroup 
$\langle \sigma \rangle$ generated by $\sigma$.  According to Theorem~\ref{AH},  
after passing to a finite cover, we may assume that $\alpha$ intersects 3                   
distinct dividing curves.  Similarly, we may assume that $\beta$ intersects 3 
distinct dividing curves.  Next, by applying Theorem~\ref{AH} to the trivial 
subgroup $H=\{e\}$ (or by using residual finiteness), there exists a finite 
cover $\widetilde M$ which does not contain $\sigma^k$.  Let $\widetilde\alpha$ 
be a component of the preimage of $\alpha$ which intersects a component 
$\widetilde A$ of the preimage of $A$ in an arc.  If $k$ is large enough, there 
is a component $\widetilde\beta$ of the preimage of $\beta$ which straddles 
the other boundary component of $\widetilde A$ and is disjoint from 
$\widetilde\alpha$.    This proves the claim.

Now, Proposition~\ref{add=dig} shows that attaching a bypass along $\widetilde 
\beta$ from the exterior of $\widetilde {M}\setminus \widetilde{\Sigma}$ 
produces a contact structure  isomorphic to the one obtained by digging out the 
bypass in $\widetilde {M}\setminus \widetilde{\Sigma}$  attached along 
$\widetilde \alpha$, and, in particular, this new contact  structure must be 
universally tight.  On the other hand, $\widetilde \beta$ is assumed to be the 
attaching arc for a bypass in $\widetilde {M}\setminus \widetilde{\Sigma}$, and 
it follows  that attaching a bypass along $\widetilde \beta$ to the outside of 
$\widetilde {M}\setminus \widetilde{\Sigma}$ produces an overtwisted disk.    
This is a contradiction. \end{proof}

If $\xi$ is a universally tight contact structure and $\Sigma$ is a component of 
$\bdry M$ which admits a full marking, then let $\Gamma(S)$ be the union of the
dividing curves of $\Sigma$ which are straddled by the (full) straddling set 
$S$.  In view of Proposition~\ref{straddled-unstraddled}, $\Gamma(S)$ is an 
invariant of $\xi$.  

We now consider the case where $M=\Sigma \times [0,1]$ has convex boundary 
which consists of $\Sigma_i=\Sigma\times\{i\}$, $i=0,1$, and where each 
$\Gamma_{\Sigma_i}$ is just a single pair of parallel dividing curves.  We 
proved the following classification theorem in \cite{HKM3}.

\begin{thm}\label{classification}
Let $\Sigma$ be a closed oriented surface of genus $g \geq 2$ and 
$M=\Sigma\times [0,1]$.  Suppose that $\Gamma_{\Sigma_i}$, $i=0,1$, is a pair 
$\gamma_i\sqcup \gamma_i'$ of parallel nonseparating curves which cobound an 
annulus $A_i \subset \Sigma_i$ and that 
$\Gamma_{\Sigma_0}\not=\Gamma_{\Sigma_1}$. Choose a characteristic foliation 
$\mathcal{F}$ on $\bdry M$ which is adapted to $\Gamma_{\Sigma_0}\sqcup 
\Gamma_{\Sigma_1}$.  Then there exist, up to isotopy rel boundary, exactly 4 
tight contact structures which satisfy the boundary condition $\mathcal{F}$, and 
all of them  are universally tight.  Moreover:

\be
\item  For each of the 4 tight contact structures, both $\Sigma_0$ and 
$\Sigma_1$ are fully straddled, and the 4 cases exactly correspond to the 
4 possible choices of $\Gamma(S)$, consisting of pairs of straddled curves, one 
from each $\Gamma_{\Sigma_i}$. 

\item Let $\delta_i$ be a closed Legendrian curve on $\Sigma_i$ which has 
geometric intersection $|\delta_i\cap \gamma_i|=1$.  Then there exists a 
(degenerate) bypass in $\Sigma\times I$ along $\delta_i$ which straddles the 
straddled curve of $\Gamma_{\Sigma_i}$.

\ee

\end{thm}

\s\n
{\it Remark.} The focus of this paper is on surfaces of genus greater than one. 
However, the definitions of straddled, unstraddled, and full marking make sense 
on a torus with only slight modification to allow for the necessity of more than 
two parallel dividing curves.  The statements and proofs above apply to tori, 
with one exception:  the uniqueness statement for universally tight contact 
structures on $T^2 \times I$ with a given full marking. The classification 
theorem in this  setting requires an additional nonnegative integer invariant 
called the {\em  torsion} (see \cite{Gi99a}, \cite{Co99a}, \cite{HKM2}).

\section{Gluing along surfaces which admit full markings}  
\label{section:gluing}

The atoroidal hypothesis in the following Gluing Theorem allows us to assume $M$ 
is hyperbolic and apply Proposition~\ref{straddled-unstraddled}.  A 
substantially similar argument can be made without the atoroidal hypothesis by 
using results of Long-Niblo~\cite{LN} instead of Allman-Hamilton~\cite{AH}.  

\begin{thm}[Gluing Theorem] \label{closedgluing} 
Let $(M,\xi)$ be an atoroidal contact manifold which is extremal along a closed, 
convex, incompressible surface $\Sigma$, and let $\Sigma^0$ and $\Sigma^1$ be 
the boundary components of $M \setminus \Sigma$ corresponding to $\Sigma$.  
Suppose that, for each $i=0,1$, $\Sigma^i$ admits a full marking in $M \setminus 
\Sigma$ with a straddling set $S^i$, and the dividing curves of $\Sigma^i$ are 
straddled if and  only if they are unstraddled in $\Sigma^{1-i}$.  Suppose 
further that $S^0 \cup S^1$ is an embedded, essential family of arcs in 
$\Sigma$. Then $(M,\xi)$ is universally tight if and only if $(M \setminus 
\Sigma, \xi|_{M \setminus \Sigma})$ is universally tight. 
\end{thm}

\begin{proof}
The proof that $(M,\xi)$ is tight also applies to finite covers of $(M, \xi)$, 
in light of Proposition~\ref{covers}.   Atoroidal Haken 3-manifolds are 
hyperbolic, and hence have residually finite fundamental groups.  Therefore, any 
overtwisted disk that exists in the universal cover also exists in some finite 
cover, and the proof below will also imply that $(M, \xi)$ is universally tight.

The general strategy for proving tightness is explained in detail in \cite{HKM2} 
and \cite{H3}, so we will only provide a brief summary.  Arguing by 
contradiction, we assume $(M,\xi)$ contains an overtwisted disk $D$.  There 
exists a sequence $\Sigma=\Sigma_0, \Sigma_1,\dots,\Sigma_n$ of isotopic 
surfaces, where each step is a single bypass attachment and $\Sigma_n$ is 
disjoint from $D$. Since we can extricate $D$ from (isotopic copies of) $\Sigma$ 
in stages, if we show that each $\xi|_{M\setminus \Sigma_i}$ is tight, this 
implies the tightness of $\xi|_M$.  Now, for universal tightness, we pass to a 
large finite cover $\widetilde M$ of $M$ and extricate some lift $\widetilde D$ 
of $D$ from the preimage $\widetilde \Sigma$ of $\Sigma$.  Lifting to a cover 
has the following advantages: 

\s\n
\be

\item A bypass which is $\#\Gamma$-increasing can be made trivial by using the 
residual finiteness of $M$.

\item A bypass whose attaching arc intersects only two distinct curves but is 
not trivial or $\#\Gamma$-increasing can be made to intersect three distinct 
curves by Theorem~\ref{AH}. 

\ee

\s\n
Therefore, by lifting as necessary so that each bypass satisfies Conditions 1 
and 2 above, we have a sequence 
$$(\widetilde M_0=M, \widetilde\Sigma_0 = \Sigma, \widetilde\xi_0 = \xi), 
(\widetilde M_1, \widetilde \Sigma_1, \widetilde \xi_1), \dots, (\widetilde M_n, 
\widetilde\Sigma_n, \widetilde \xi_n),$$ 
where $\widetilde\Sigma_i\subset \widetilde M_i$, $\widetilde M_{i+1}$ is a 
finite cover of $\widetilde M_i$,  $\widetilde \xi_i$ is the pullback of $\xi$ 
to $\widetilde M_i$, $\widetilde \Sigma_{i+1}$ is obtained from the preimage of 
$\widetilde \Sigma_i$ by attaching a bypass along an arc which straddles three 
distinct curves,  and a lift $\widetilde D$ of $D$ in $\widetilde M_n$ is  
disjoint from $\widetilde \Sigma_n$.  If we can show that 
$\widetilde\xi_i|_{\widetilde M_i\setminus \widetilde \Sigma_i}$ is universally 
tight for all $i$, we are done.  This involves making sure that at each step the 
contact structure admits a full marking and satisfies the conditions of the 
theorem --- therefore we need to exercise extra care when choosing the covers.

Consider the attachment of the first bypass.  We start with $\Sigma$ and $M$, 
and take a finite cover which satisfies Conditions 1 and 2.  (To simplify 
notation, we will still write $M$ for the finite cover of $M$ and $\Sigma$ for 
the preimage of $\Sigma$.)  Then the bypass $\mathcal{B}_\alpha$ has an 
attaching arc $\alpha$ which intersects three distinct dividing curves of 
$\Sigma$.  Denote the copies of $\Sigma$ in $M\setminus \Sigma$ by $\Sigma^i$, 
$i=0,1$, and let $S^i$ be the straddling set for $\Sigma^i$.  By 
Proposition~\ref{covers}, $\Sigma$ admits a full marking on both sides (i.e., on 
$\Sigma^0$ and on $\Sigma^1$), and it is clear that the union $S^0\cup S^1$ is 
still an embedded essential family of arcs.  After passing to a larger finite 
cover and removing duplicate stradding arcs which straddle the same dividing 
curve as in Proposition~\ref{covers}:

\s\n
\be
\item[3.]  The arc of attachment of a bypass can be made disjoint from the 
straddling sets $S^i$, $i=0,1$, by using residual finiteness as in the proof of 
Proposition~\ref{straddled-unstraddled}. 
\ee
\s\n

Suppose we dig out $\mathcal{B}_\alpha$ from the $\Sigma^0$ side and reattach 
$\mathcal{B}_\alpha$ along $\alpha$ to the $\Sigma^1$ side.  This gives us 
$M\setminus \Sigma'$, where $\Sigma'$ is a surface parallel to and disjoint 
from $\Sigma$.   By Proposition~\ref{straddled-unstraddled}, the curve $\gamma$ 
straddled by $\alpha$ must be in $\Gamma(S^0)$ .  Since $\gamma$ is straddled by 
an element $\beta_0$ of $S^0$ if and only if $\gamma$ is unstraddled by an 
element of $S^1$, it follows that the parallel curve $\gamma'$ is straddled by 
an arc $\beta_1 \in S^1$.

It is clear that digging $\mathcal{B}_\alpha$ preserves universal tightness. 
Since $\alpha$ and $\beta_1$ are disjoint (by Condition~3),  
Proposition~\ref{add=dig} tells us that, on $M\setminus \Sigma$, digging 
$\mathcal{B}_\alpha$ out is equivalent to attaching a bypass 
$\mathcal{B}_{\beta_1}'$ along $\beta_1$ onto the $\Sigma^0$ side, i.e., 
$M\setminus (\Sigma\cup \mathcal{B}_\alpha)$ is contactomorphic to 
$M'=\overline{(M\setminus \Sigma)\cup N(\mathcal{B}'_{\beta_1})}$, where $N(F)$ 
is a small neighborhood of $F$.  We will write $\bdry M'=(\Sigma')^0\sqcup 
\Sigma^1$.  By Proposition~\ref{bypass-keeps-straddled}, since $S^0 \cup 
\{\beta_1\}$ is essential, $(\Sigma')^0\subset M'$ admits a full marking 
obtained by dropping $\beta_0$ from $S^0$.  On the other hand, by 
Proposition~\ref{add=dig} again, attaching a bypass $\mathcal{B}_\alpha'$ to 
$\Sigma^1\subset M'$ along $\alpha$ is equivalent to digging a bypass 
$\mathcal{B}_{\beta_1}$ along $\beta_1$.  In other words, $M'\cup 
N(\mathcal{B}_\alpha')$ is contactomorphic to $M'\setminus 
\mathcal{B}_{\beta_1}$ and also to $M\setminus \Sigma'$. Therefore, 
$\xi|_{M\setminus \Sigma'}$ is universally tight and admits full markings $S^0 - 
\{\beta_0 \}$ on the $(\Sigma')^0$ side and $S^1 - \{\beta_1 \}$ on the 
$(\Sigma')^1$ side (by applying Proposition~\ref{bypass-keeps-straddled} again).  
The union of the two straddling sets is an embedded essential family of arcs. 

We can now inductively construct  $(\widetilde M_i, \widetilde \Sigma_i, \widetilde
\xi_i)$ that admit full markings on both sides which satisfy the conditions of 
the theorem, by choosing finite covers where Conditions 1, 2, and 3 are met. 
This ensures universal tightness of each step and finishes the proof.  
\end{proof}

\section{The Gabai-Eliashberg-Thurston Theorem}  \label{section:get}

\begin{proof} [Proof of Theorem~\ref{GET-closed}]
According to Gabai \cite{Ga}, there is a well-groomed sutured manifold 
decomposition of $M$, 
$$M\stackrel{\Sigma}{\sa} (M_1,\gamma_1)\stackrel{S_1}{\sa}\cdots 
\stackrel{S_{n-1}}{\sa}(M_n,\gamma_n)=\cup (B^3,S^1 \times I),$$ 
where $\Sigma$ is a nonseparating surface.  Since $M$ and $\Sigma$ are closed, 
$\gamma_1=\emptyset$.  For $i \ge 1$, $S_i$ may be chosen to have nonempty 
boundary (see \cite{HKM2}, Theorem~1.3 for a statement of this version of 
Gabai's theorem).  It follows that for $i \ge 2$, $(M_i, \gamma_i)$ has annular 
sutures, that is, all sutures are annuli and every component of $\partial M_i$ 
contains at least one suture.

By the results of \cite{HKM2}, the sutured manifold decomposition
$$(M_2,\gamma_2)\stackrel{S_2} {\sa}(M_3,\gamma_3) \stackrel{S_{3}}{\sa} 
\cdots\stackrel{S_{n-1}}{\sa}(M_n,\gamma_n)=\cup (B^3,S^1 \times I),$$ 
gives rise to a convex decomposition 
$$(M_2,\Gamma_2)\stackrel{(S'_2,\sigma_2)}{\sa} (M_3,\Gamma_3) 
\stackrel{(S'_{3},\sigma_{3})}{\sa} \cdots\stackrel{(S'_{n-1},\sigma_{n-1})} 
{\sa} (M_n,\Gamma_n)=\cup (B^3,S^1),$$ 
and then (and this requires that each of the $(M_i, \Gamma_i)$ above have 
annular sutures) $(M_2,\Gamma_2)$ carries a universally tight contact structure 
by Theorem~6.1 of \cite{HKM1}.  Note that the proof of Theorem~6.1 in 
\cite{HKM1} uses the perturbation of a taut foliation into a universally tight 
contact structure;  we give a foliation-theory-free proof of the same fact in 
\cite{HKM2}.

Gabai's construction gives $(M_1,\gamma_1=\emptyset)\stackrel{S_1} 
{\sa}(M_2,\gamma_2)$, so to apply Theorem~6.1 of \cite{HKM1}, it is necessary to 
produce a convex structure $(M_1, \Gamma_1)$ with annular sutures (in 
particular, $\Gamma_1\neq \emptyset$) and a splitting surface $(S'_1,\sigma_1)$ 
such that $(M_1,\Gamma_1)\stackrel{(S'_1,\sigma_1)}{\sa}(M_2,\Gamma_2)$. 
Let $\Sigma^0$ and $\Sigma^1$ be the components of $\bdry M_1$ corresponding 
to the original splitting surface $\Sigma$.  Since $S_1$ is well-groomed, the 
components of $S_1 \cap \Sigma^i$, $i=0,1$, are parallel oriented nonseparating 
curves in the isotopy class $s_i$.  Let $A_i$ be  an annular neighborhood of 
a curve dual to $s_i$, and denote $\partial A_i = \delta_i \sqcup \delta_i'$. 
Define $\Gamma_1 = \delta_0 \sqcup \delta_0' \sqcup \delta_1 \sqcup \delta_1'$. 
The convex structure $(M_1, \Gamma_1)$ is defined by decreeing that $\Sigma^i 
\setminus A_i \subset (\Sigma^i)_+$ and $A_i \subset (\Sigma^i)_-$ if and only 
if the orientation induced from $\Sigma \setminus A$ agrees with the outward 
pointing normal orientation on $\Sigma^i \setminus A_i$.

Now, $(S'_1,\sigma_1)$ is defined so that $S_1'=S_1$ and $\sigma_1$ is the 
unique dividing set which is $\bdry$-parallel and gives rise to $(M_2,\Gamma_2)$ 
after the splitting.   The bypasses corresponding to the $\bdry$-parallel 
dividing curves straddle curves of $\Gamma_1$.   Due to the well-grooming of 
$S_1$, there is a unique choice of straddled curve for each pair 
$\delta_i\sqcup \delta_i'$.  It now follows that there is a universally tight 
contact structure on $(M_1, \Gamma_1)$ and that $\bdry M_1$ admits a full 
marking with $\Gamma(S)= \delta_0\sqcup\delta_1$, for example.

Next consider $M= M_1\cup (\Sigma\times[0,1])$, where we identify $\Sigma^i$ 
with $\Sigma\times\{i\}$, $i=0,1$.  By Theorem~\ref{classification}, there is a 
(unique) universally tight contact structure on $\Sigma \times I$ with 
$\Gamma_{\bdry (\Sigma\times I)}=\delta_0 \sqcup \delta_0' \sqcup \delta_1 
\sqcup \delta_1'$ and $\Gamma(S) = \delta_0'\sqcup \delta_1'$.   We apply 
Theorem~\ref{closedgluing} twice to obtain a universally tight contact structure 
on $M$, as follows.   First consider the gluing of $\Sigma^0$ to 
$\Sigma\times\{0\}$.   Denote the straddling set for $\Sigma^0$ by $S^0$ and the 
straddling set for $\Sigma\times\{0\}$ by $S^1$.  Then $S^0$ consists of one 
closed curve isotopic to $s_0$ (note that the corresponding bypass is 
degenerate).  Choose $S^1$ so that it consists of one closed curve dual 
to $\delta_0$ and $S^0\cup S^1$ is an embedded, essential family on 
$\Sigma^0=\Sigma\times\{0\}$.  This is easily arranged since the genus of 
$\Sigma^0$ is $\geq 2$.    Finally, by (2) of Theorem~\ref{classification}, 
there exists a (degenerate) bypass in $\Sigma\times I$ along the unique closed 
curve in $S^1$ which straddles $\delta_0'$.  Now we can apply 
Theorem~\ref{closedgluing}.  The second gluing is identical. \end{proof}

\end{document}